\documentclass{article}
\usepackage{amssymb}
\usepackage{amsfonts}
\usepackage{amsmath}

\input{tcilatex}
\begin{document}

\bigskip

\begin{titlepage}
\title{\bf A New Type Warped Product Metric in Contact Geometry}
\author{Ahmet MOLLAOGULLARI {\footnote{E-mail: ahmet\_m@comu.edu.tr}} \\
  {\small Department of Matehematics, Faculty of Art and Science,   \c{C}anakkale Onsekiz Mart University,}
\\{\small  17100  \c{C}anakkale, Turkey.}
\\ \c{C}etin CAMCI \footnote{  E-mail : ccamci@comu.edu.tr} \\
 {\small Department of Matehematics, Faculty of Art and Science, 
 \c{C}anakkale Onsekiz Mart University, }
\\{\small  17100 \c{C}anakkale, Turkey.}
}

\maketitle
\begin{abstract}

In this study, we show that there is an $\alpha $-Sasakian structure
on product manifold $M_{1}\times \beta (I)$ where $M_{1}$ is a \ Kaehlerian
manifold that has exact 1-form and $\beta (I)$ is an open curve. After then, we define a new type warped product metric
to study the warped product of almost Hermitian manifolds with almost contact metric manifolds.

\end{abstract}

\end{titlepage}

\section{Introduction}

Contact geometry is a theoretical subject which has so many applications in
the fields of science such as physics and engineering. From thermodynamics
to optics, from electrics to motion equations, it has an important place in
many areas \cite{Gieges}.Many studies have been carried out on contact
manifolds and contact geometry which became increasingly important in the
20th century. In addition, the symplectic geometry, Hermitian manifolds and
Kaehlerian manifolds which have serious applications in many fields, are
also important topics in mathematics. That's why, to consider them together,
we studied the product of contact and complex manifolds in this stdudy.

The product of contact and complex manifolds have been an interesting area
for mathematicians. Pandey (1981) studied the necessary and sufficient
condition for a product manifold to be almost complex, Kaehler and almost
Tachibana manifold. Caprusi (1984) showed that the product of two almost
contact manifolds can be an almost Hermitian manifold. In addition, he
determined the necessary and sufficient conditions for a product manifold to
be a Hermitian, Kaehler, almost Kaehler and nearly-Kaehler. Zegga et al.
(2017) investigated the product of Kaehler and Sasaki manifolds via the
metric

\begin{equation}
g=g_{1}+g_{2}+\theta \otimes \theta +\theta \otimes \eta +\eta \otimes \theta
\end{equation}

and obtained a Sasaki structure on this product manifold. Also Gherici at
al. (2019) showed that the product of real numbers and a Kaehler manifold
with a Kaehler form has a Sasakian structure via the metric

\begin{equation}
g(X,Y)=h(X,Y)+\theta (X)\theta (Y)\text{ , \ }g(X,\partial _{r})=\theta (X)%
\text{, \ }g(\partial _{r},\partial _{r})=1
\end{equation}

Also there are valuable studies about warped product which takes an
important place in our paper in \cite{arslan},\cite{atceken},\cite{Chen} and 
\cite{cihan}

Because metric is fundemental tool to study a product manifols, we defined
the following metric on the product of almost Hermitian manifold $%
(M_{1},J,g_{1})$ which has an exact Kahlerian form and the open curve $\beta
:I\rightarrow E^{n}$ to get contact and Sasakian structures on $%
M=M_{1}\times \beta (I)$ with the structure $(\bar{\phi},\bar{\xi},\bar{\eta}%
):$

\begin{equation}
g=g_{1}+\bar{\eta}\otimes \bar{\eta}
\end{equation}

After them, for an almost Hermitian manifold ($M_{1},J,g_{1})$ with exact
1-form $\omega $ which satisfies $d\omega =\Phi _{1}$ and for an almost
contact metric manifold $(M_{2},\phi ,\xi ,\eta ,g_{2}),$ we defined a new
type warped product metric on the structure $\left( M_{1}\times M_{2},\bar{%
\phi},\bar{\xi},\bar{\eta}\right) $ as following:

\begin{equation}
g=g_{1}+f^{2}\left( g_{2}-\eta \otimes \eta -\bar{\eta}\otimes \bar{\eta}%
\right)
\end{equation}

By (4), we studied the connection on $M_{1}\times M_{2}$ when $M_{1}$ is an
almost Hermitian manifold with an axact 1-form and $M_{2}$ is almost cantact
metric manifold, metric manifold and K-contact manifold respectiveley. \
Finally, we determined the necessary and sufficient conditions for a curve
to be a geodesic in the warped product manifold $M=M_{1}\times M_{2}$

\section{Preliminaries}

\bigskip

\textbf{Definition 2.1.}\cite{B} Let $M$ be a $(2n+1)-$dimensional manifold
and $\phi ,\xi ,\eta $ \ be\ $(1,1),(1,0)$ and $(0,1)$ tensor fields
respectively on $M$. If these tensor fields satisfy the equatins (3) and (4)
then ($\phi ,\xi ,\eta )$ is called an almost contact structure on $M$ and $%
(M,\phi ,\xi ,\eta )$ is called an almost contact manifold.%
\begin{eqnarray}
\phi ^{2}X &=&-X+\eta (X)\xi \\
\eta (\xi ) &=&1
\end{eqnarray}%
From this definition it can be deduced that $\phi (\xi )=0$ and $\eta \circ
\phi =0.$

\bigskip

\bigskip \textbf{Definition 2.2.} \textit{\cite{B1}}Assume that $(M,\phi
,\xi ,\eta )$ is an almost contact manifold with dimension $2n+1$. If $g$ is
a Riemannian or Lorentzian metric and$\ $ satisfies the equation 
\begin{equation}
g(\phi (X),\phi (Y))=g(X,Y)-\eta (X)\eta (Y),\text{ }\forall X,Y\in \chi (M)
\end{equation}%
\qquad\ then $(\phi ,\xi ,\eta ,g)$ is called an almost contact metric
structure and $(M,\phi ,\xi ,\eta ,g)$ is called an almost contact metric
manifold.

\bigskip

\textbf{Definiton 2.3.} \cite{B1}Assume that $(M,\phi ,\xi ,\eta ,g)$ is a $%
(2n+1)-$dimensional almost contact metric manifold. If the equation%
\begin{equation}
d\eta (X,Y)=g(X,\phi (Y))
\end{equation}%
is satisfied on $M$ then $(M,\phi ,\xi ,\eta ,g)$ is called a contact metric
manifold.

\bigskip

\textbf{Lemma 2.1.}\cite{B} Let $(M,\phi ,\xi ,\eta ,g)$ be an almost
contact metric manifold then for all $\ X,Y\in \chi (M)$ the followings hold:

\textit{i)} $g(X,\xi )=\eta (X)$ \ ,

\textit{ii)}$g(\phi (X),\xi )=0$ ,

\textit{iii)}$g(\phi (X),Y)+g(\phi (Y),X)=0$.

\bigskip

\textbf{Definition 2.4.} \cite{YK}Let $(M,\phi ,\xi ,\eta ,g)$ be an almost
contact metric manifold. The fundemantal 2-form $\Phi $ of this manifold is
defined as following: 
\begin{equation}
\Phi :\chi (M)\times \chi (M)\rightarrow C^{\infty }(M,%
%TCIMACRO{\U{211d} }%
%BeginExpansion
\mathbb{R}
%EndExpansion
),\text{ \ \ \ }\Phi \left( X,Y\right) =g(X,\phi Y).
\end{equation}

\bigskip

\textbf{Definition 2.5. }\cite{YK}Let $(M,J,g)$ be an almost complex
manifold. The fundemental 2-form of this manifold is defined as following:%
\begin{equation}
\Phi :\chi (M)\times \chi (M)\rightarrow C^{\infty }(M,%
%TCIMACRO{\U{211d} }%
%BeginExpansion
\mathbb{R}
%EndExpansion
),\ \ \Phi \left( X,Y\right) =g(X,JY).
\end{equation}

\bigskip

\textbf{Definition 2.6.} Assume that $(M,\phi ,\xi ,\eta ,g)$ is an almost
contact metric manifold. M is called a Sasakian manifold if the structure $%
(\phi ,\xi ,\eta ,g)$ satisfies the following equation:%
\begin{equation}
(\nabla _{X}\phi )Y=g(X,Y)\xi -\eta (Y)X,\text{ \ \ }\forall X,Y\in \chi (M).
\end{equation}%
More generally an $\alpha -$Sasakian structure may be defined by the
requirement%
\begin{equation}
(\nabla _{X}\phi )Y=\alpha (g(X,Y)\xi -\eta (Y)X),\text{ \ \ }\forall X,Y\in
\chi (M)
\end{equation}%
where $\alpha $ is a non-zero constant.

\bigskip

\textbf{Theorem 2.1.} \cite{YK}Let $(M,\phi ,\xi ,\eta )$ be an almost
contact manifold then for all $X,Y\in \chi (M)$ the following equation holds:%
\begin{equation}
2d\eta (X,Y)=X\eta (Y)-Y\eta (X)-\eta \lbrack X,Y]
\end{equation}%
where [,] stands for Lie bracket.

\bigskip

\textbf{Definition 2.7.} \cite{Pan}On a product manifold such as $%
M=M_{1}\times M_{2},$ the Lie derivative of two vector fields $X=\left(
X_{1},X_{2}\right) ,Y=\left( Y_{1},Y_{2}\right) \in $ $\chi (M),$ is defined
as following: 
\begin{equation}
\lbrack \text{ },]:\chi (M)\times \chi (M)\rightarrow \chi (M),\text{ \ \ }%
[X,Y]=\left( [X_{1},Y_{1}],[X_{2},Y_{2}]\right)
\end{equation}

\bigskip

\textbf{Definition 2.8.} \cite{Pan} Let $f_{1}:M_{1}\rightarrow 
%TCIMACRO{\U{211d} }%
%BeginExpansion
\mathbb{R}
%EndExpansion
$ and $f_{2}:M_{2}\rightarrow 
%TCIMACRO{\U{211d} }%
%BeginExpansion
\mathbb{R}
%EndExpansion
$ be differantiable functions on the manifolds $M_{1}$ and $M_{2}$
respectiveley. For $X_{1}\in \chi (M_{1})$ and $X_{2}\in \chi (M_{2})$ the
derivative of $f_{1}+f_{2}$ and $f_{1}.f_{2}$ in the direction of the vector
field $X=\left( X_{1},X_{2}\right) \in \chi (M_{1}\times M_{1})$ are defined
as follows:\ 
\begin{equation}
X[f_{1}+f_{2}]=X_{1}[f_{1}]+X_{2}[f_{2}]\text{ }
\end{equation}%
\begin{equation}
X[f_{1}.f_{2}]=f_{2}X_{1}[f_{1}]+f_{1}X_{2}[f_{2}]
\end{equation}

\bigskip

\textbf{Theorem 2.2.} \cite{B1}Let (M,g) be a Riemannian manifold and $%
\nabla $ be the Levi-Civita connection. Then for all \ $X,Y,Z\in \chi (M)$
the following equaition holds:%
\begin{equation}
\begin{array}{c}
2g(\nabla _{X}Y,Z)=Xg(Y,Z)+Yg(X,Z)-Zg(X,Y)-g(X,[Y,Z]) \\ 
+g(Y,[Z,X])+g(Z,[X,Y])%
\end{array}%
\end{equation}%
This equiation is also called as Koszul Formula.

\bigskip

\section{\textbf{Contact and Sasakian Structures on }$M=M_{1}\times \protect%
\beta (I)$\textbf{\ }}

In this section we examine the conditions on $M_{1}$ and the open curve $%
\beta :I\rightarrow E^{n}$ to get contact and Sasakian structures on $%
M=M_{1}\times \beta (I)$. \ We suppose that ($M_{1},J,g_{1})$ is an almost
Hermitian manifold of dimension $2n$ which has an exact 1-form $\omega $
which satisfies%
\begin{equation}
d\omega (X_{1},Y_{1})=\Phi _{1}(X_{1},Y_{1})=g_{1}(X_{1},JY_{1})\text{ , \ }%
\forall X_{1},Y_{1}\in \chi (M_{1}).
\end{equation}%
Therefore, there exists a $\left(
x_{1},x_{2},...,x_{n},y_{1},y_{2},...,y_{n}\right) $ coordinate system for
any neighborhood U of $M_{1}$ such that%
\begin{eqnarray}
J\partial y_{i} &=&\partial x_{i} \\
J\partial x_{i} &=&-\partial y_{i}
\end{eqnarray}%
equations hold for $i=1,2,...,n$ where $\frac{\partial }{\partial y_{i}}%
=\partial y_{i}$ and $\frac{\partial }{\partial x_{i}}=\partial x_{i}$ . We
can choose the set 
\begin{equation*}
\left\{ 2\partial x_{1},2\partial x_{2},...,2\partial x_{n},2\partial
y_{1},2\partial y_{2},..,2\partial y_{n}\right\}
\end{equation*}%
as an orthonormal base for $M_{1}$ without lose of generality. Therefore the
set$\{\frac{1}{2}dx_{1},\frac{1}{2}dx_{2},...,\frac{1}{2}dx_{n},\frac{1}{2}%
dy_{1},\frac{1}{2}dy_{2},...,\frac{1}{2}dy_{n}\}$is the dual of the chosen
base. For the next computations we take the dual base as $\{\frac{1}{2}%
\omega _{1},\frac{1}{2}\omega _{2},...,\frac{1}{2}\omega _{n},\frac{1}{2}%
\omega _{n+1},\frac{1}{2}\omega _{n+2},...,\frac{1}{2}\omega _{2n}\}$ So $%
\omega $ can be written in the form of dual base vectors such as:

\begin{equation}
\omega =\frac{1}{2}f_{1}\omega _{1}+\frac{1}{2}f_{2}\omega _{2}+...+\frac{1}{%
2}f_{2n}\omega _{2n}=\frac{1}{2}\overset{2n}{\underset{i=1}{\sum }}%
f_{i}\omega _{i}
\end{equation}%
where $f_{i}$ $\in C^{\infty }(M_{1},%
%TCIMACRO{\U{211d} }%
%BeginExpansion
\mathbb{R}
%EndExpansion
).$ Thus we have the following equations:%
\begin{equation}
\omega (2\partial x_{i})=f_{i}\ \ ,\ \ \omega (2\partial y_{i})=f_{n+i}\ \ \
i=1,2,...,n
\end{equation}

Furthermore because $d\omega =\Phi _{1}$ the following equatins hold for the
coefficient functions of $\omega $ :%
\begin{eqnarray}
\frac{\partial f_{j}}{\partial x_{i}}-\frac{\partial f_{i}}{\partial x_{j}}
&=&0,\text{ \ }i\neq j \\
\frac{\partial f_{n+i}}{\partial y_{j}}-\frac{\partial f_{n+j}}{\partial
y_{i}} &=&0,\text{ \ }i\neq j \\
\frac{\partial f_{n+j}}{\partial x_{i}}-\frac{\partial f_{i}}{\partial y_{j}}
&=&0,\text{ \ }i\neq j \\
\frac{\partial f_{n+i}}{\partial x_{i}}-\frac{\partial f_{i}}{\partial y_{i}}
&=&1,\text{ \ }i=j
\end{eqnarray}

\textbf{Theorem 3.1. }The triple $\left( \bar{\phi},\bar{\xi},\bar{\eta}%
\right) $ is an almost contact sructure on $M$ where%
\begin{equation}
\bar{\phi}:\chi (M)\rightarrow \chi (M),\ \bar{\phi}\left( X_{1},hT\right)
=\left( JX_{1},-\alpha \omega (JX_{1})T\right)
\end{equation}%
\begin{equation}
\bar{\eta}:\chi (M)\rightarrow C^{\infty }(M,%
%TCIMACRO{\U{211d} }%
%BeginExpansion
\mathbb{R}
%EndExpansion
),\bar{\eta}\left( X_{1},hT\right) =\alpha \omega (X_{1})+h
\end{equation}%
\begin{equation}
\bar{\xi}=\left( 0,T\right)
\end{equation}%
Here T stands for the unit tangent vector field of $\beta $ and $h$ is a
function as $h:I\rightarrow 
%TCIMACRO{\U{211d} }%
%BeginExpansion
\mathbb{R}
%EndExpansion
$

\bigskip

\textbf{Proof \ }For $X=(X_{1},hT)$ and $\bar{\xi}\in \chi (M)$ we have the
following equations:%
\begin{eqnarray}
\bar{\phi}^{2}(X) &=&(-X_{1},\alpha \omega (X_{1})T)  \notag \\
&=&-X+\bar{\eta}(X)\bar{\xi}
\end{eqnarray}%
and\qquad 
\begin{equation}
\bar{\eta}(\bar{\xi})=\bar{\eta}(0,T)=1
\end{equation}%
By means of definiton 2.1., $\left( M,\bar{\phi},\bar{\xi},\overset{\mathbf{-%
}}{\eta }\right) $ is an almost contact manifold.

\bigskip

\textbf{Theorem 3.2. }$\left( M,\bar{\phi},\bar{\xi},\bar{\eta},g\right) $
is an almost contact manifold with the metric 
\begin{eqnarray}
g &:&\chi (M)\times \chi (M)\rightarrow C^{\infty }(M,%
%TCIMACRO{\U{211d} }%
%BeginExpansion
\mathbb{R}
%EndExpansion
)  \notag \\
\text{\ }g\left( \left( X_{1},h_{1}T\right) ,\left( Y_{1},h_{2}T\right)
\right) &=&g_{1}(X_{1},Y_{1})+\bar{\eta}\left( X_{1},h_{1}T\right) \bar{\eta}%
\left( Y_{1},h_{2}T\right)
\end{eqnarray}

\textbf{Proof }For $X=\left( X_{1},h_{1}T\right) \in \chi (M)$ and $Y=\left(
Y_{1},h_{2}T\right) \in \chi (M)$ we obtain 
\begin{eqnarray*}
g(\bar{\phi}(X),\bar{\phi}(Y)) &=&g_{1}(JX_{1},JY_{1})+\bar{\eta}(\bar{\phi}%
X).\bar{\eta}(\bar{\phi}Y) \\
&=&g_{1}(X_{1},Y_{1}) \\
&=&g(X,Y)-\bar{\eta}(X)\bar{\eta}(Y)
\end{eqnarray*}%
From definition 2.2., it's concluded that $\left( M,\bar{\phi},\bar{\xi},%
\overset{\mathbf{-}}{\eta },g\right) $ is an almost contact metric manifold.

\bigskip

\textbf{Lemma 3.1. }For all $X=(X_{1},h_{1}T)$, $Y=(Y_{1},h_{2}T)\in \chi
(M) $ there is a relationship between the fundemental 2-forms of $\left( M,%
\bar{\phi},\bar{\xi},\bar{\eta},g\right) $ and ($M_{1},J,g_{1})$ as 
\begin{equation}
\bar{\Phi}\left( X,Y\right) =\Phi _{1}(X_{1},Y_{1}).
\end{equation}

\bigskip

\textbf{Proof }One can easily see by using the deifnitions of fundemental
2-forms of $\left( M_{1},J,g_{1})\right) $ and $\left( M,\bar{\phi},\bar{\xi}%
,\bar{\eta},g\right) $

\bigskip

\textbf{Theorem 3.3.} $\left( M,\bar{\phi},\bar{\xi},\bar{\eta},g\right) $
is a contact metric manifold if and only if $\alpha =1.$

\bigskip

\textbf{Proof }For $X=\left( X_{1},h_{1}T\right) \in \chi (M)$ and $Y=\left(
Y_{1},h_{2}T\right) \in \chi (M)$ we get,%
\begin{eqnarray*}
2d\bar{\eta}(X,Y) &=&X\bar{\eta}(Y)-Y\bar{\eta}(X)-\bar{\eta}[X,Y] \\
&=&X_{1}\left[ \alpha \omega Y_{1}\right] -Y_{1}\left[ \alpha \omega X_{1}%
\right] -\omega \lbrack X_{1},Y_{1}] \\
&=&2ad\omega (X_{1},Y_{1}) \\
&=&2a\Phi _{1}(X_{1},Y_{1}) \\
&=&2a\Phi (X,Y)
\end{eqnarray*}%
So $\bar{\eta}$ is exact 1-form if and only if $\alpha =1.$ Therefore $%
\left( M,\bar{\phi},\bar{\xi},\bar{\eta},g\right) $ is a contact metric
manifold if and only if $\alpha =1.$

\bigskip

Now, we determine the conditions for the product manifold $M$ to be a
Sasakian manifold. For this aim we consider the base of $M$ as $\left\{ 
\overset{\mathbf{-}}{\partial x_{1}},\overset{\mathbf{-}}{\partial x_{2}}%
,...,\overset{\mathbf{-}}{\partial x_{n}},\overset{\mathbf{-}}{\partial y_{1}%
},\overset{\mathbf{-}}{\partial y_{2}},..,\overset{\mathbf{-}}{\partial
y_{n},}\overset{\mathbf{-}}{\partial z}\right\} $ where%
\begin{equation*}
\overset{\mathbf{-}}{\partial x_{i}}=(\partial x_{i},0)\text{ ,\ \ \ \ \ }%
\overset{\mathbf{-}}{\partial y_{i}}=(\partial y_{i},0)\text{ \ ,\ \ \ \ \ }%
\overset{\mathbf{-}}{\partial z}=\frac{1}{2}\left( 0,T\right)
\end{equation*}%
By Theorem 3.2. we can write the followings:%
\begin{eqnarray*}
g(2\partial x_{i},\overset{\mathbf{-}}{2\partial x_{j}}) &=&g_{1}(2\partial
x_{i},2\partial x_{j})+\bar{\eta}\left( 2\partial x_{i},0\right) \bar{\eta}%
\left( 2\partial x_{j},0\right) \\
&=&\delta _{ij}+[\alpha \omega (2\partial x_{i})][\alpha \omega (2\partial
x_{i})] \\
&=&\delta _{ij}+\alpha ^{2}f_{i}f_{j}
\end{eqnarray*}%
Similarly others are found as%
\begin{eqnarray*}
g(\overset{\mathbf{-}}{2\partial x_{i}},\overset{\mathbf{-}}{2\partial y_{j}}%
) &=&\alpha ^{2}f_{i}f_{n+j} \\
g(\overset{\mathbf{-}}{2\partial x_{i}},\overset{\mathbf{-}}{\partial z})
&=&\alpha f_{i} \\
g(2\partial y_{i},2\partial y_{j}) &=&\delta _{ij}+\alpha ^{2}f_{n+i}f_{n+j}
\\
g(\overset{\mathbf{-}}{2\partial y_{i}},2\partial z\text{\ }) &=&\alpha
f_{n+i} \\
g(2\partial z,2\partial z) &=&1
\end{eqnarray*}%
So the matrix corresponding to the metric $g$ on the base $\left\{ \overset{%
\mathbf{-}}{\partial x_{1}},\overset{\mathbf{-}}{\partial x_{2}},...,\overset%
{\mathbf{-}}{\partial x_{n}},\overset{\mathbf{-}}{\partial y_{1}},\overset{%
\mathbf{-}}{\partial y_{2}},..,\overset{\mathbf{-}}{\partial y_{n},}\overset{%
\mathbf{-}}{\partial z}\right\} $ is founds as%
\begin{equation*}
\left[ \mathbf{g}_{ij}\right] \mathbf{=}\frac{1}{4}\left[ 
\begin{array}{ccc}
\mathbf{\delta }_{ij}\mathbf{+\alpha }^{2}\mathbf{f}_{i}\mathbf{f}_{j} & 
\mathbf{\alpha }^{2}\mathbf{f}_{i}\mathbf{f}_{n+j} & \mathbf{\alpha f}_{i}
\\ 
\mathbf{\alpha }^{2}\mathbf{f}_{j}\mathbf{f}_{n+i} & \mathbf{\delta }_{ij}%
\mathbf{+\alpha }^{2}\mathbf{f}_{n+i}\mathbf{f}_{n+j} & \mathbf{\alpha f}%
_{n+i} \\ 
\mathbf{\alpha f}_{j} & \mathbf{\alpha f}_{n+j} & \mathbf{1}%
\end{array}%
\right] _{(2n+1)\times (2n+1)}
\end{equation*}%
And inverse of $\left[ \mathbf{g}_{ij}\right] $ ;%
\begin{equation*}
\left[ g_{ij}\right] ^{-1}=\left[ g^{ij}\right] =4\left[ 
\begin{array}{ccc}
\mathbf{\delta }_{ij} & \mathbf{0} & \mathbf{-\alpha f}_{i} \\ 
0 & \mathbf{\delta }_{ij} & -\mathbf{\alpha f}_{n+j} \\ 
\mathbf{-\alpha f}_{i} & -\mathbf{\alpha f}_{n+j} & 1+\mathbf{\alpha }^{2}%
\overset{2n}{\underset{i=1}{\sum }}\mathbf{f}_{i}^{2}%
\end{array}%
\right]
\end{equation*}%
For the later, we have to calculate the Christoffel symbols for this metric.
It is known that 2$^{nd}$ kind of Christoffel symbols can be calculated by
the formula \cite{Blair2} 
\begin{equation*}
\Gamma _{ij}^{k}=\frac{1}{2}\overset{2n+1}{\underset{h=1}{\sum }}%
g^{kh}\left( g_{_{ih,j}}+g_{_{hj,i}}-g_{_{ij,h}}\right)
\end{equation*}%
So non-zore of them are calculated as followings:%
\begin{eqnarray*}
\Gamma _{_{ij}}^{^{^{n+i}}} &=&\frac{1}{2}\alpha ^{2}f_{j}\text{ },\text{ \ }%
\Gamma _{ij}^{^{n+j}}=\frac{1}{2}\alpha ^{2}f_{i}\text{ },\text{ \ \ }\Gamma
_{ii}^{^{n+i}}=\alpha ^{2}f_{i} \\
\Gamma _{_{ij}}^{^{2n+1}} &=&\alpha \frac{\partial f_{j}}{\partial x_{i}}-%
\frac{1}{2}\alpha ^{3}\left( f_{i}f_{n+j}+f_{j}f_{n+i}\right) \text{ \ } \\
\Gamma _{_{i(n+j)}}^{^{j}} &=&-\frac{1}{2}\alpha ^{2}f_{i}\text{, \ }\Gamma
_{i(n+j)}^{^{n+i}}=\frac{1}{2}\alpha ^{2}f_{n+j} \\
\Gamma _{i(n+j)}^{^{2n+1}} &=&\alpha \frac{\partial f_{i}}{\partial y_{j}}+%
\frac{1}{2}\alpha \delta _{ij}+\frac{1}{2}\alpha ^{3}\left(
f_{i}f_{j}-f_{n+i}f_{n+j}\right) \\
\Gamma _{i(2n+1)}^{^{n+i}} &=&\frac{1}{2}\alpha \text{ \ , }\Gamma
_{i(2n+1)}^{^{2n+1}}=-\frac{1}{2}\alpha ^{2}f_{n+i}\text{ \ } \\
\text{\ }\Gamma _{(n+i)(n+j)}^{^{i}} &=&-\frac{1}{2}\alpha ^{2}f_{n+j}\text{
, \ }\Gamma _{(n+i)(n+j)}^{^{j}}=-\frac{1}{2}\alpha ^{2}f_{n+i} \\
\Gamma _{(n+i)(n+j)}^{^{2n+1}} &=&\alpha \frac{\partial f_{n+j}}{\partial
y_{i}}+\frac{1}{2}\alpha ^{3}\left( f_{i}f_{n+j}+f_{j}f_{n+i}\right) \\
\Gamma _{(n+i)(2n+1)}^{^{i}} &=&-\frac{1}{2}\alpha \text{ \ \ , \ \ }\Gamma
_{(n+i)(2n+1)}^{^{2n+1}}=\frac{1}{2}\alpha ^{2}f_{i}
\end{eqnarray*}%
If we define the vector fileds $e_{i}=2\partial \overset{-}{y_{i}}-\alpha
f_{n+i}\overset{-}{\xi }$\ on $M,$ we get $\bar{\phi}(e_{i})=2\partial 
\overset{-}{x_{i}}-\alpha f_{i}\overset{\mathbf{-}}{\xi }$, $\ \bar{\eta}%
(e_{i})=0$ and also $\bar{\eta}(\bar{\phi}(e_{i}))=0.$Thus we obtain%
\begin{eqnarray*}
g\left( e_{i},e_{j}\right) &=&g_{1}(2\partial y_{i},2\partial y_{j})+\bar{%
\eta}(e_{i})\bar{\eta}(e_{j}) \\
&=&\delta _{ij}
\end{eqnarray*}%
Similarly we can obtain the following equations:%
\begin{eqnarray}
g\left( \bar{\phi}(e_{i}),\bar{\phi}(e_{j})\right) &=&\delta _{ij} \\
g\left( e_{i},\bar{\phi}(e_{j})\right) &=&0 \\
g\left( e_{i},\bar{\xi}\right) &=&0 \\
g\left( \bar{\phi}(e_{i}),\bar{\xi}\right) &=&0
\end{eqnarray}%
Thefore the set $\{e_{1},e_{2},...,e_{n},\bar{\phi}(e_{1}),\bar{\phi}%
(e_{2}),...,\bar{\phi}(e_{n}),\bar{\xi}$\ \ $\}$ is an orthonormal base for $%
M.$

\bigskip

\textbf{Lemma 3.2.} For the set $\{e_{1},e_{2},...,e_{n},\bar{\phi}(e_{1}),%
\bar{\phi}(e_{2}),...,\bar{\phi}(e_{n}),\bar{\xi}$\ $\},$we have the
followings:%
\begin{eqnarray*}
\nabla _{e_{i}}\text{ }e_{j} &=&\nabla _{\bar{\phi}(e_{i})}\text{ }\bar{\phi}%
(e_{j})=\nabla _{\bar{\xi}\text{\ }}\bar{\xi}\text{\ }=0 \\
\nabla _{e_{i}}\text{ }\bar{\phi}(e_{j}) &=&\alpha \delta _{ij}\bar{\xi}%
\text{\ }=-\nabla _{\bar{\phi}(e_{i})}e_{j} \\
\nabla _{\bar{\xi}\text{\ }}\text{ }e_{i} &=&-\alpha \bar{\phi}%
(e_{i})=\nabla _{e_{i}}\bar{\xi}\text{\ } \\
\nabla _{\bar{\xi}\text{\ }}\text{ }\bar{\phi}(e_{i}) &=&\alpha e_{i}=\nabla
_{\bar{\phi}(e_{i})}\bar{\xi}\text{\ }
\end{eqnarray*}

\textbf{Proof }For $e_{i}=2\partial \overset{-}{y_{i}}-\alpha f_{n+i}\overset%
{\mathbf{-}}{\xi },$ we have

\begin{eqnarray*}
\nabla _{\bar{\xi}\text{\ }}e_{i} &=&\dsum\limits_{k=1}^{2n+1}\Gamma
_{_{(2n+1)(n+i)}}^{^{k}}\partial k-a\left( \overset{\mathbf{-}}{\xi }%
[f_{n+i}]\overset{\mathbf{-}}{\xi }+f_{n+i}\nabla _{\overset{\mathbf{-}}{\xi 
}}\overset{\mathbf{-}}{\xi }\right) \\
&=&-a(2\partial \overset{-}{x_{i}}-\alpha f_{i}\overset{\mathbf{-}}{\xi }) \\
&=&-a\bar{\phi}(e_{i}).
\end{eqnarray*}%
Others can be shown similarly.

\bigskip

\textbf{Theorem} \textbf{3.4.} \bigskip Let ($M_{1},J,g_{1})$ be a Kahlerian
manifold with the exact 1-form. Then $\left( M,\bar{\phi},\bar{\xi},\bar{\eta%
},g\right) $ is a Sasakian manifold if and only if $\alpha =1$

\textbf{Proof} \ Let $\tilde{\nabla}$ be the Levi-Civita connetction on $%
M_{1}.$ Using Lemma 3.3 for $X=\left( X_{1},h_{1}T\right) \in \chi (M)$ and $%
Y=\left( Y_{1},h_{2}T\right) \in \chi (M)$ one can find the Levi-Civita
connetction on $M$ as following:

\begin{equation}
{\small \nabla }_{X}{\small Y=}\left( 
\begin{array}{c}
\tilde{\nabla}_{X_{1}}Y_{1}-a\bar{\eta}(Y)JX_{1}-a\bar{\eta}(X)JY_{1}, \\ 
\left( X\bar{\eta}(Y)-a\omega (\tilde{\nabla}%
_{X_{1}}Y_{1})+g_{1}(JX_{1},Y_{1})+a^{2}\bar{\eta}(Y)\omega (JX_{1})+a^{2}%
\bar{\eta}(X)\omega (JY_{1})\right) T%
\end{array}%
\right) .
\end{equation}%
So $\nabla _{X}\left( \bar{\phi}Y\right) $ is calculatd as

\begin{eqnarray*}
\nabla _{X}\left( \bar{\phi}Y\right) &=&\left( 
\begin{array}{c}
\tilde{\nabla}_{X_{1}}JY_{1}-a\bar{\eta}(\bar{\phi}Y)JX_{1}-a\bar{\eta}%
(X)J^{2}Y_{1}, \\ 
\left( X\bar{\eta}(\bar{\phi}Y)-a\omega (\tilde{\nabla}%
_{X_{1}}JY_{1})+g_{1}(JX_{1},JY_{1})+a^{2}\bar{\eta}(\bar{\phi}Y)\omega
(JX_{1})+a^{2}\bar{\eta}(X)\omega (J^{2}Y_{1})\right) T%
\end{array}%
\right) . \\
&=&\left( \tilde{\nabla}_{X_{1}}JY_{1}+a\bar{\eta}(X)Y_{1},\left( -a\omega (%
\tilde{\nabla}_{X_{1}}JY_{1})+g_{1}(X_{1},Y_{1})-a^{2}\bar{\eta}(X)\omega
(Y_{1})\right) T\right)
\end{eqnarray*}%
Similarly $\bar{\phi}(\nabla _{X}Y)$ \ can be abotained as

\begin{equation}
\bar{\phi}(\nabla _{X}Y)=\left( 
\begin{array}{c}
J\tilde{\nabla}_{X_{1}}Y_{1}+a\bar{\eta}(Y)X_{1}+a\bar{\eta}(X)Y_{1} \\ 
\left( -\omega \left( J\tilde{\nabla}_{X_{1}}Y_{1}\right) -a\bar{\eta}%
(Y)\omega X_{1}-a\bar{\eta}(X)\omega Y_{1}\right) T%
\end{array}%
\right)
\end{equation}%
So we get

\begin{eqnarray}
\left( \nabla _{X}\bar{\phi}\right) Y &=&\nabla _{X}\left( \bar{\phi}%
Y\right) -\bar{\phi}(\nabla _{X}Y)  \notag \\
&=&\left( 
\begin{array}{c}
\left( \tilde{\nabla}_{X_{1}}J\right) Y_{1}-a\bar{\eta}(X)Y_{1} \\ 
a(-\omega \left( \tilde{\nabla}_{X_{1}}J\right) Y_{1}+g_{1}(X_{1},Y_{1})+%
\bar{\eta}(Y)\omega X_{1})T%
\end{array}%
\right) .
\end{eqnarray}%
Because $M_{1}$ is Kahlerian $\tilde{\nabla}J=0$ and

\begin{equation}
g_{1}(X_{1},Y_{1})=g(X,Y)-a^{2}\omega (X_{1})\omega (Y_{1})-a\omega
(X_{1})h_{2}-a\omega (Y_{1})h_{1}-h_{1}h_{2}
\end{equation}%
then equatiaon (38) turns to

\begin{eqnarray}
\left( \nabla _{X}\bar{\phi}\right) Y &=&a\left( \bar{\eta}(X)Y_{1},\left(
g_{1}(X_{1},Y_{1})+\bar{\eta}(Y)\omega X_{1}\right) T\right)  \notag \\
&=&a(g(X,Y)\bar{\xi}-\bar{\eta}(Y)X)+(-a^{2}a\omega (Y_{1})h_{1}+\bar{\eta}%
(Y)h_{1}-ah_{1}h_{2})\bar{\xi}.
\end{eqnarray}%
By forward calculations it is seen that

\begin{equation}
-a^{2}a\omega (Y_{1})h_{1}+\bar{\eta}(Y)h_{1}-ah_{1}h_{2}=0
\end{equation}%
So $\left( \nabla _{X}\bar{\phi}\right) Y$ is found as

\begin{equation}
\left( \nabla _{X}\bar{\phi}\right) Y=a(g(X,Y)\bar{\xi}-\bar{\eta}(Y)X).
\end{equation}%
Equation (42) tells that $M=M_{1}\times \beta (I)$ is an $a-$Sasakian
manifold with the given structure. And it is a Sasakian manifold if and only
if $\alpha =1$

$\bigskip $

\textbf{Example} Consider the Kaehlerian manifold $\left( 
%TCIMACRO{\U{211d} }%
%BeginExpansion
\mathbb{R}
%EndExpansion
^{4},J,g_{1}\right) $ where the metric is given as%
\begin{equation}
g_{1}=\frac{1}{2}dx_{1}^{2}+\frac{1}{2}dx_{2}^{2}+\frac{1}{2}dy_{1}^{2}+%
\frac{1}{2}dy_{2}^{2}
\end{equation}%
\ and the almost complex structure is given as 
\begin{equation}
J(\partial x_{i})=\partial y_{i}\ ,\ J(\partial y_{i})=-\partial x_{i}i=1,2.
\end{equation}%
Taking $\omega =-\frac{1}{2}x_{1}dy_{1}-\frac{1}{2}x_{2}dy_{2}$ results d$%
\omega =\frac{1}{2}\left( dy_{1}\wedge dx_{1}+dy_{2}\wedge dx_{2}\right) .$
One can see that d$\omega =\Phi _{1}.$ We consider the curve $\beta :(0,2\pi
)\rightarrow E^{2},$ $\beta (t)=(\cos t,\sin t)$ with tangent $T=(-\sin
t,\cos t)=2\frac{d}{dt}.$ So for $M=%
%TCIMACRO{\U{211d} }%
%BeginExpansion
\mathbb{R}
%EndExpansion
^{4}\times \beta (I),$ the structure $\left( M,\bar{\phi},\bar{\xi},\bar{\eta%
},g\right) $ is an a-Sasakian manifold with the followings:%
\begin{equation}
\bar{\phi}:\chi (M)\rightarrow \chi (M),\bar{\phi}\left( X_{1},2h\frac{d}{dt}%
\right) =\left( JX_{1},-2a\omega (JX_{1})\frac{d}{dt}\right)
\end{equation}

\begin{equation}
\bar{\eta}:\chi (M)\rightarrow C^{\infty }(M,%
%TCIMACRO{\U{211d} }%
%BeginExpansion
\mathbb{R}
%EndExpansion
),\text{ }\bar{\eta}\left( X_{1},2h\frac{d}{dt}\right) =a\omega (X_{1})+h
\end{equation}

\begin{equation}
g=\frac{1}{4}\left[ 
\begin{array}{ccccc}
1 & 0 & 0 & 0 & 0 \\ 
0 & 1 & 0 & 0 & 0 \\ 
0 & 0 & 1+a^{2}x_{1}^{2} & a^{2}x_{1}x_{2} & -ax_{1} \\ 
0 & 0 & a^{2}x_{1}x_{2} & 1+a^{2}x_{2}^{2} & -ax_{2} \\ 
0 & 0 & -ax_{1} & -ax_{2} & 1%
\end{array}%
\right]
\end{equation}

\section{A New Type Warped Product Metric}

\bigskip In this section, we define a new type warped product metric and
examine some product manifolds by using this metric.

\textbf{Theorem 4.1. } Let ($M_{1},J,g_{1})$ be an almost Hermitian manifold
with exact 1-form $\omega $ which satisfies $d\omega =\Phi _{1}$ and $%
(M_{2},\phi ,\xi ,\eta ,g_{2})$ be an almost contact metric manifold. Then $%
\left( M_{1}\times M_{2},\bar{\phi},\bar{\xi},\bar{\eta},g\right) $ is an
almost contact metric manifold with the following tensor feilds;

\begin{eqnarray}
\bar{\phi}(X_{1},X_{2}) &=&\left( JX_{1},\phi X_{2}-a\omega (JX_{1})\xi
\right) , \\
\bar{\eta}\left( X_{1},X_{2}\right) &=&a\omega (X_{1})+\eta \left(
X_{2}\right) , \\
\bar{\xi} &=&\left( 0,\xi \right) \\
g\left( (X_{1},X_{2}),\left( Y_{1},Y_{2}\right) \right)
&=&g_{1}(X_{1},X_{2})+f^{2}\left( g_{2}(X_{2},Y_{2})-\eta (X_{2})\eta
(Y_{2})-\bar{\eta}(X)\bar{\eta}(Y)\right)
\end{eqnarray}%
where $X=\left( X_{1},X_{2}\right) \in \chi (M),Y=\left( Y_{1},Y_{2}\right)
\in \chi (M)$ , $a\in 
%TCIMACRO{\U{211d} }%
%BeginExpansion
\mathbb{R}
%EndExpansion
$ and $f\in C^{\infty }(M_{2},%
%TCIMACRO{\U{211d} }%
%BeginExpansion
\mathbb{R}
%EndExpansion
)$ is a positive function.\bigskip

\textbf{Proof }By forward calculation it can be found that 
\begin{eqnarray*}
\bar{\phi}^{2}(X_{1},X_{2}) &=&\left( J^{2}X_{1},\phi ^{2}(X_{2}\right)
-a\omega (JX_{1})\phi (\xi )-a\omega (J^{2}\left( X_{1}\right) \xi ) \\
&=&\left( -X_{1},-X_{2}+\eta (X_{2})\xi +a\omega (\left( X_{1}\right) \xi
\right) \\
&=&-X+\bar{\eta}(X)\bar{\xi}
\end{eqnarray*}

\begin{equation*}
\bar{\eta}(\bar{\xi})=\bar{\eta}(0,\xi )=1
\end{equation*}%
Then it's understood that $\left( M_{1}\times M_{2},\bar{\phi},\bar{\xi},%
\bar{\eta}\right) $ is an almost contact manifold. Also, because 
\begin{eqnarray*}
g\left( \bar{\phi}(X),\bar{\phi}(Y)\right) &=&g\left( \left( JX_{1},\phi
X_{2}-a\omega (JX_{1})\xi \right) ,\left( JY_{1},\phi Y_{2}-a\omega
(JY_{1})\xi \right) \right) \\
&=&g_{1}(X_{1},Y_{1})+f^{2}g_{2}(\phi (X_{2}),\phi (Y_{2})) \\
&=&g_{1}(X_{1},Y_{1})+f^{2}(g_{2}(X_{2},Y_{2})-\eta (X_{2})\eta (Y_{2})) \\
&=&g(X,Y)-\bar{\eta}(X)\bar{\eta}(Y)
\end{eqnarray*}

then $\left( M_{1}\times M_{2},\bar{\phi},\bar{\xi},\bar{\eta},g\right) $ is
an almost contact metric manifold.

\bigskip

\textbf{Theorem 4.2. }Let ($M_{1},J,g_{1})$ be an almost Hermitian manifold
with exact 1-form $\omega $ which satisfies $d\omega =\Phi _{1}$ and $%
(M_{2},\phi ,\xi ,\eta ,g_{2})$ be an almost contact metric manifold. The
following equations hold for all $X,Y\in \chi (M_{1})$ and $U,V,W\in \chi
(M_{2})$ and for the Levi-Civita the Levi-Civita connections $\bar{\nabla},%
\tilde{\nabla}$ and $\nabla $ on $M_{1},M_{2}$ and $M=M_{1}\times M_{2}$
respectively:{}

\begin{enumerate}
\item[i.] $\nabla _{(X,0)}(Y,0)=\left( \bar{\nabla}_{X}Y-\alpha ^{2}\omega
(X)J(Y)-\alpha ^{2}\omega (Y)J(X),\alpha \left( -\omega (\bar{\nabla}%
_{X}Y)+\alpha ^{2}\omega (X)\omega (JY)+\alpha ^{2}\omega (Y)\omega (JX)+%
\frac{\lambda }{2})\xi \right) \right) $

where $\lambda =\lambda (X,Y)=X\omega (Y)+Y\omega (X)+\omega ([X,Y])$

\item[ii.] $g\left( \nabla _{(X,0)}(0,U),(Y,W)\right) =fX[f]g_{2}\bigskip
\left( \phi (U),\phi (W)+2a\eta (U)d\omega (X,Y)+2a\omega (X)d\eta
(U,W)\right) $

\item[iii.] $g\bigskip \left( \nabla _{(0,U)}(0,V),(Y,W)\right)
=f^{2}g_{2}\left( \tilde{\nabla}_{U}V,W\right) +(1-f^{2})\eta (V)d\eta
(U,W)+(1-f^{2})\eta (U)d\eta (V,W)+\frac{\theta }{2}\bar{\eta}(Y,W)-\frac{%
\theta }{2}f^{2}\eta (W)-fY[f]g_{2}(\phi (U),\phi (V))$\newline
where $\theta =\theta (U,V)=U\eta (V)+V\eta (U)+\eta \lbrack U,V]$
\end{enumerate}

\bigskip

\textbf{Proof \ i) }Using Koszul formula for $(Z,W)\in \chi (M)$ we get

\begin{eqnarray*}
2g\left( \nabla _{(X,0)}(Y,0),(Z,W)\right) &=&2g_{1}\left( \bar{\nabla}%
_{X}Y,Z\right) +2\alpha ^{2}\omega (X)d\omega (Y,Z)+2\alpha ^{2}\omega
(Y)d\omega (X,Z) \\
&&+\alpha \left( \alpha \omega (Z)+\eta (W)\right) \left( X\omega
(Y)+Y\omega (X)+\omega ([X,Y])\right) \\
&=&2g_{1}\left( \bar{\nabla}_{X}Y,Z\right) -2\alpha ^{2}\omega
(X)g_{1}(JY,Z)-2\alpha ^{2}\omega (Y)g_{1}(JX,Z)+\alpha \bar{\eta}%
(Z,W)\lambda \\
&=&2g_{1}(\bar{\nabla}_{X}Y-\alpha ^{2}\omega (X)JY-\alpha ^{2}\omega
(Y)JX,Z)+\alpha \lambda g((0,\xi ),(Z,W))
\end{eqnarray*}%
On the other by putting $\bar{\nabla}_{X}Y-\alpha ^{2}\omega (X)JY-\alpha
^{2}\omega (Y)JX=\mu $ one can see that

\begin{equation*}
2g_{1}(\bar{\nabla}_{X}Y-\alpha ^{2}\omega (X)JY-\alpha ^{2}\omega
(Y)JX,Z)=2g\left( (\mu ,-\alpha \omega (\mu )\xi ,(Z,W)\right)
\end{equation*}%
So we get

\begin{eqnarray*}
g\left( \nabla _{(X,0)}(Y,0),(Z,W)\right) &=&g\left( (\mu ,-\alpha \omega
(\mu )\xi ,(Z,W)\right) +\frac{\alpha }{2}\lambda g((0,\xi ),(Z,W)) \\
&=&g\left( \left( \mu ,-\alpha \omega (\mu )\xi \right) +\frac{\alpha }{2}%
\lambda (0,\xi ),\left( Z,W\right) \right)
\end{eqnarray*}%
Because the metric is positive definite it's concluded that

\begin{eqnarray*}
\nabla _{(X,0)}(Y,0) &=&\left( \mu ,-\alpha \omega (\mu )\xi \right) +\frac{%
\alpha }{2}\lambda (0,\xi ) \\
&=&\left( \mu ,\alpha \left( -\alpha \omega (\mu )+\frac{\lambda }{2}\right)
\xi \right) \\
&=&\left( \bar{\nabla}_{X}Y-\alpha ^{2}\omega (X)J(Y)-\alpha ^{2}\omega
(Y)J(X),\alpha \left( -\omega (\bar{\nabla}_{X}Y)+\alpha ^{2}\omega
(X)\omega (JY)+\alpha ^{2}\omega (Y)\omega (JX)+\frac{\lambda }{2}\right)
\xi \right)
\end{eqnarray*}%
\textbf{ii)} \textbf{\ }Using Koszul formula for $(Y,W)\in \chi (M)$ we get

\begin{eqnarray*}
2g\left( \nabla _{(X,0)}(0,U),(Y,W)\right)
&=&(X,0)g((0,U),(Y,W))+(0,U)g((X,0),(Y,W))+(Y,W)g((X,0),(0,U)) \\
&&-g((X,0),[(0,U),(Y,W)]+g((0,U),[(Y,W),(0,U)]+g((Y,W),[(X,0),(0,U)]) \\
&=&2fX[f]g_{2}(\phi (U),\phi (W))+\alpha \eta (U)X\omega (Y)+\alpha \omega
(X)U\eta (W)-\alpha \eta (U)Y\omega (X) \\
&&-\alpha \omega (X)W\eta (U)-\alpha \omega (X)\eta \lbrack U,W]+\alpha \eta
(U)\omega \lbrack Y,X] \\
&=&2fX[f]g_{2}\bigskip \left( \phi (U),\phi (W)+2a\eta (U)d\omega
(X,Y)+2a\omega (X)d\eta (U,W)\right)
\end{eqnarray*}%
Then it's found that

\begin{equation*}
g\left( \nabla _{(X,0)}(0,U),(Y,W)\right) =fX[f]g_{2}\bigskip \left( \phi
(U),\phi (W))+a\eta (U)d\omega (X,Y)+a\omega (X)d\eta (U,W)\right)
\end{equation*}%
\textbf{iii)} \textbf{\ }Again using Koszul formula for $(Y,W)\in \chi (M)$
we get

\begin{eqnarray*}
2g\bigskip \left( \nabla _{(0,U)}(0,V),(Y,W)\right)
&=&(0,U)g((0,V),(Y,W))+(0,V)g((0,U),(Y,W))-(Y,W)g((0,U),(0,V)) \\
&&-g((0,U),(0,[V,W]))+g((0,V),(0,[W,U]))+g((Y,W),(0,[U,V])) \\
&=&f^{2}g_{2}\left( \tilde{\nabla}_{U}V,W\right) +2(1-f^{2})\eta (V)d\eta
(U,W)+2(1-f^{2})\eta (U)d\eta (V,W) \\
&&+(1-f^{2})\theta \eta (W)+\alpha \omega (X)\theta -2Yf[Y]g_{2}\left( \phi
(U),\phi (V)\right)
\end{eqnarray*}%
Because $\bar{\eta}(Y,W)=\alpha \omega (X)+\eta (W)$ it's calculated that

\begin{eqnarray*}
g\bigskip \left( \nabla _{(0,U)}(0,V),(Y,W)\right) &=&f^{2}g_{2}\left( 
\tilde{\nabla}_{U}V,W\right) +(1-f^{2})\eta (V)d\eta (U,W)+(1-f^{2})\eta
(U)d\eta (V,W) \\
&&+\frac{\theta }{2}\bar{\eta}(Y,W)-\frac{f^{2}}{2}\theta \eta
(W)-fY[f]g_{2}(\phi (U),\phi (V))\newline
\end{eqnarray*}

\textbf{Theorem 4.3. }Let ($M_{1},J,g_{1})$ be an almost Hermitian manifold
with exact 1-form $\omega $ which satisfies $d\omega =\Phi _{1}$ and $%
(M_{2},\phi ,\xi ,\eta ,g_{2})$ be a contact metric manifold. For all $X\in
\chi (M_{1})$ and $U\in \chi (M_{2})$ the following equation holds:

\begin{eqnarray}
\nabla _{(X,0)}(0,U) &=&\nabla _{(0,U)}(X,0) \\
&=&\left( -\alpha \eta (U)JX,\frac{X[f]}{f}\phi ^{2}(U)-\frac{\alpha }{f^{2}}%
\omega (X)\phi (U)+\alpha ^{2}\eta (U)\omega (JX)\xi \right)
\end{eqnarray}

\textbf{Proof }If $(M_{2},\phi ,\xi ,\eta ,g_{2})$ is a contact metric
manifold then it's also an almost contact manifold. Then (ii) of Theorem
4.2. holds. Besides, becasue $(M_{2},\phi ,\xi ,\eta ,g_{2})$ is a contact
metric manifold, then for all $U,W\in \chi (M_{2})$

\begin{equation}
d\eta (U,W)=g_{2}(U,\phi (W)).
\end{equation}%
Putting $d\eta (U,W)=g_{2}(U,\phi (W))$ and $d\omega (X,Y)=g_{1}(X,JY)$ in
(ii)

\begin{equation}
g\left( \nabla _{(X,0)}(0,U),(Y,W)\right) =fX[f]g_{2}\bigskip (\phi (U),\phi
(W))+a\eta (U)g_{1}(X,JY)+a\omega (X)g_{2}(U,\phi (W)).
\end{equation}%
On the other hand using the definiton of the metric $g$ we find

\begin{equation}
g(\bar{\phi}(X,0),(Y,W))=-g_{1}(X,JY)
\end{equation}

\begin{equation}
g((0,U),\bar{\phi}(Y,W))=f^{2}g_{2}(U,\phi (W)).
\end{equation}%
Because g is self adjoint we write

\begin{eqnarray}
g((0,U),\bar{\phi}(Y,W)) &=&-g(\bar{\phi}(0,U),(Y,W))  \notag \\
&=&-g((0,\phi U),(Y,W))  \notag \\
&=&-f^{2}g_{2}(U,\phi (W))
\end{eqnarray}%
Replacing $\phi U$ instead $U$ it's seen that

\begin{equation}
g_{2}(\phi (U),\phi (W))=-\frac{1}{f^{2}}g((0,\phi ^{2}(U)),(Y,W)).
\end{equation}%
Applying (56), (57) and (58) in (55) we conclude

\begin{equation*}
g\left( \nabla _{(X,0)}(0,U),(Y,W)\right) =g\left( \left( -\frac{X[f]}{f}%
(0,\phi ^{2}(U)-\frac{a}{f^{2}}\omega (X)(0,\phi (U))-a\eta (U)\bar{\phi}%
(X,0)\right) ,(Y,W)\right)
\end{equation*}%
Since tihs equation is provided for all $(Y,W)\in \chi (M)$ then\bigskip 
\begin{equation*}
\nabla _{(X,0)}(0,U)=\left( -\alpha \eta (U)JX,\frac{X[f]}{f}\phi ^{2}(U)-%
\frac{\alpha }{f^{2}}\omega (X)\phi (U)+\alpha ^{2}\eta (U)\omega (JX)\xi
\right)
\end{equation*}%
Besides, using $[(X,0),[0,U])=(0,0)$ and $[(X,0),[0,U])=\nabla
_{(X,0)}(0,U)-\nabla _{(0,U)}(X,0)$ we get

\begin{equation*}
\nabla _{(X,0)}(0,U)=\nabla _{(0,U)}(X,0)
\end{equation*}

\textbf{Theorem 4.4. }Let ($M_{1},J,g_{1})$ be an almost Hermitian manifold
with exact 1-form $\omega $ which satisfies $d\omega =\Phi _{1}$ and $%
(M_{2},\phi ,\xi ,\eta ,g_{2})$ be a K-contact manifold. The following
equation holds for all $U,V\in \chi (M_{2})$ and $U,V,W\in \chi (M_{2})$ and
for the Levi-Civita connections $\bar{\nabla},\tilde{\nabla}$ and $\nabla $
on $M_{1},M_{2}$ and $M=M_{1}\times M_{2}$ respectively:

\begin{equation}
\nabla _{(0,U)}(0,V)=\left( 
\begin{array}{c}
-fg_{2}(\phi (U),\phi (V))\func{grad}f, \\ 
\tilde{\nabla}_{U}V+\frac{f^{2}-1}{f^{2}}\left( \eta (V)\phi (U)+\eta
(U)\phi (V)+\alpha fg_{2}\left( \phi (U),\phi (V)\right) \omega (\func{grad}%
f)\xi \right)%
\end{array}%
\right)
\end{equation}

\textbf{Proof }Since $(M_{2},\phi ,\xi ,\eta ,g_{2})$ is a K-contact
manifold it's also a contact metric manifold and satisfies

\begin{eqnarray}
d\eta (U,W) &=&g_{2}(U,\phi (W)) \\
d\eta (V,W) &=&g_{2}(V,\phi (W))
\end{eqnarray}%
$.$

Because the preconditions in (iii) of Theorem 4.2. are provided, using (61)
and (622) in (iii), it's found that;

\begin{eqnarray}
2g\bigskip \left( \nabla _{(0,U)}(0,V),(Y,W)\right) &=&g_{2}\left( 2f^{2}%
\tilde{\nabla}_{U}V+2(f^{2}-1)\eta (V)\phi (U)+2(f^{2}-1)\eta (U)\phi
(V)-f^{2}\theta \xi ,W\right) \\
&&-2fg_{1}(\func{grad}f,Y)g_{2}(\phi (U),\phi (V))+\theta \bar{\eta}(Y,W) 
\notag
\end{eqnarray}%
In this equatin if we put $2f^{2}\tilde{\nabla}_{U}V+2(f^{2}-1)\eta (V)\phi
(U)+2(f^{2}-1)\eta (U)\phi (V)-f^{2}\theta \xi =\mu $ then it's seen that

\begin{equation}
g_{2}\left( \mu ,W\right) =\frac{1}{f^{2}}g\left( \left( 0,\mu \right)
,\left( Y,W\right) \right)
\end{equation}%
On the other hand it's seen that

\begin{equation}
g(\left( \func{grad}f,0\right) ,\left( Y,W\right) )=g_{1}(\func{grad}%
f,Y)+\alpha \omega \left( \func{grad}f\right) g((0,\xi ),(Y,W))
\end{equation}%
So we find

\begin{equation}
Y[f]=g_{1}(\func{grad}f,Y)=g((\func{grad}f,-\alpha \omega \left( \func{grad}%
f\right) \xi ),((Y,W)))
\end{equation}%
Writing $\bar{\eta}(Y,W)=g\left( (0,\xi ),(Y,W)\right) $ and applying
equations (64),(65) and (66) in (63) we conclude that

\begin{equation*}
\nabla _{(0,U)}(0,V)=\left( 
\begin{array}{c}
-fg_{2}(\phi (U),\phi (V))\func{grad}f, \\ 
\tilde{\nabla}_{U}V+\frac{f^{2}-1}{f^{2}}\left( \eta (V)\phi (U)+\eta
(U)\phi (V)+\alpha fg_{2}\left( \phi (U),\phi (V)\right) \omega (\func{grad}%
f)\xi \right)%
\end{array}%
\right)
\end{equation*}

\bigskip \textbf{Theorem 4.5. }Let ($M_{1},J,g_{1})$ be an almost Hermitian
manifold with exact 1-form $\omega $ which satisfies $d\omega =\Phi _{1}$
and $(M_{2},\phi ,\xi ,\eta ,g_{2})$ be a K-contact manifold. Let $\gamma
:I\rightarrow M_{1}$ and $\beta :I\rightarrow M_{2}$ geodesic curves on $%
M_{1}$ and $M_{2}$ respetively. \ For $X=\gamma ^{\prime }$ and $V=\beta
^{\prime }$ vector fileds, the curve $\left( \gamma ,\beta \right)
:I\rightarrow M_{1}\times M_{2}$ is a geodesic on $M_{1}\times M_{2}$ if and
only if

\begin{enumerate}
\item[i.] \ $2\alpha \bar{\eta}(X,V)JX+fg_{2}(\phi (V),\phi (V))\func{grad}%
f=0$

\item[ii.] $-\frac{2X[f]}{f}\phi ^{2}V+2\left( -\frac{1}{f^{2}}\bar{\eta}%
(X,V)+\eta (V)\right) \phi (V)$

$+\alpha \left( 2\alpha \bar{\eta}(X,V)+X\omega (JX)+X\omega (X)+fg_{2}(\phi
(V),\phi (V))\omega (\func{grad}f)\right) \xi =0$
\end{enumerate}

\textbf{Proof }Let $\bar{\nabla},\tilde{\nabla}$ and $\nabla $ be the
Levi-Civita connections on $M_{1},M_{2}$ and $M=M_{1}\times M_{2}$
respectively. Because $\gamma :I\rightarrow M_{1}$ and $\beta :I\rightarrow
M_{2}$ geodesic curves, for $X=\gamma ^{\prime }$ and $V=\beta ^{\prime }$
vector fields $\bar{\nabla}_{X}X=0$ and $\tilde{\nabla}_{V}V=0$. \ For $%
(X,V)=(\gamma ^{\prime },\beta ^{\prime })=\left( \gamma ,\beta \right)
^{\prime }$

\begin{eqnarray*}
\nabla _{(X,V)}(X,V) &=&\nabla _{(X,0)}(X,0)+\nabla _{(X,0)}(0,V)+\nabla
_{(0,V)}(X,0)+\nabla _{(0,V)}(0,V) \\
&=&\nabla _{(X,0)}(X,0)+2\nabla _{(X,0)}(0,V)+\nabla _{(0,V)}(0,V) \\
&=&\left( \bar{\nabla}_{X}X-\alpha ^{2}\omega (X)J(X)-\alpha ^{2}\omega
(X)J(X),\alpha \left( -\omega (\bar{\nabla}_{X}X)+\alpha ^{2}\omega
(X)\omega (JX)+\alpha ^{2}\omega (X)\omega (JX)+\frac{\lambda }{2})\xi
\right) \right) \\
&&+2\left( -\alpha \eta (V)JX,\frac{X[f]}{f}\phi ^{2}(V)-\frac{\alpha }{f^{2}%
}\omega (X)\phi (V)+\alpha ^{2}\eta (V)\omega (JX)\xi \right) \\
&&+\left( 
\begin{array}{c}
-fg_{2}(\phi (U),\phi (V))\func{grad}f, \\ 
\tilde{\nabla}_{V}V+\frac{f^{2}-1}{f^{2}}\left( \eta (V)\phi (V)+\eta
(V)\phi (V)+\alpha fg_{2}\left( \phi (V),\phi (V)\right) \omega (\func{grad}%
f)\xi \right)%
\end{array}%
\right) \\
&=&\left( 
\begin{array}{c}
-2\alpha \bar{\eta}(X,V)JX-fg_{2}(\phi (V),\phi (V))\func{grad}f, \\ 
-\frac{2X[f]}{f}\phi ^{2}V+\left( -\frac{2}{f^{2}}\bar{\eta}(X,V)+2\eta
(V)\right) \phi (V) \\ 
+\left( 2\alpha ^{2}\bar{\eta}(X,V)+\alpha X\omega (JX)+\alpha X\omega
(X)+\alpha fg_{2}(\phi (V),\phi (V))\omega (\func{grad}f)\right) \xi%
\end{array}%
\right)
\end{eqnarray*}

\bigskip

So to be $\nabla _{(X,V)}(X,V)=0$ , (i) and (ii) are the necessary and
sufficient conditons for the curve \bigskip $\left( \gamma ,\beta \right) $
to be a geodesic curve in $M=M_{1}\times M_{2}$

\end{document}